\documentclass{amsart}

\usepackage[english]{babel}
\usepackage{amssymb}
\usepackage{amsmath}
\usepackage{mathtools}
\usepackage[cmtip,all]{xy}
\usepackage{tikz}
\usepackage{tikz-cd}
\usepackage{enumitem}
\usepackage{changepage}
\usepackage{pdflscape}
\usepackage{geometry}
\usepackage{color}
\definecolor{linkcolor}{rgb}{0,0,0.6}
\usepackage[square,sort,comma,numbers]{natbib}
\usepackage[colorlinks=true,
pdfstartview=FitV,
linkcolor= linkcolor,
citecolor= linkcolor,
urlcolor= linkcolor,
hyperindex=true,
hyperfigures=false]
{hyperref}
\usepackage{scalerel}
\usepackage{amsthm}

\newtheorem{theorem}{Theorem}[section]
\newtheorem*{theorem*}{Theorem}
\newtheorem{lemma}[theorem]{Lemma}
\newtheorem{proposition}[theorem]{Proposition}
\newtheorem{corollary}[theorem]{Corollary}
\theoremstyle{definition}

\theoremstyle{remark}
\newtheorem{remark}[theorem]{Remark}

\numberwithin{equation}{section}

\begin{document}

\title[Nearby cycles on the ramified $\mathrm{GU}(1,n-1)$ local model]{Nearby cycles on the local model for the $\mathrm{GU}(n-1,1)$ PEL Shimura variety over a ramified prime}

\author[J. Muller]{Joseph Muller}
\address{Graduate School of Mathematical Sciences, the University of Tokyo}
\curraddr{NCTS, National Taiwan University}
\email{muller@ncts.ntu.edu.tw}

%    \subjclass is required.

%    Abstract is required.
\begin{abstract}
In this paper, we compute the cohomology sheaves of the $\ell$-adic nearby cycles on the local model of the PEL $\mathrm{GU}(n-1,1)$ Shimura variety over a ramified prime, with level given by the stabilizer of a self-dual lattice. This local model is known to have isolated singularities. If $n=2$ it has semi-stable reduction, and if $n\geq 3$ the blow-up at the singular point has semi-stable reduction. We compute the nearby cycles on the blow-up, then use proper base change to describe them on the original local model. As a result, we prove that the nearby cycles are trivial when $n$ is odd, and that only a single higher cohomology sheaf does not vanish when $n$ is even. In this case, we also describe the Galois action by computing the associated Frobenius eigenvalue.
\end{abstract}

\maketitle

\section{Introduction} 
The $\mathrm{GU}(n-1,1)$ PEL Shimura varieties can be defined over an imaginary quadratic number field $\mathbb E$ as a moduli space of abelian varieties equipped with additional structures, ie. a polarization, an $\mathbb E$-action and a level structure. Over an unramified odd prime $p$, flat integral models with parahoric level at $p$ can be defined naturally by adapting directly the moduli problem to the integral setting, see \cite{RZ} for the definition and \cite{gortz} for the proof of flatness. However, the corresponding construction over ramified primes results in an integral model which is not flat. To remedy this issue, Pappas suggested in \cite{pappas} to modify the moduli problem by adding a so-called ``wedge condition'', resulting in an integral model which, when the level at $p$ is given by the stabilizer of a self-dual lattice, is flat. We mention that the wedge condition is not enough to guarantee flatness for more general lattice stabilizer levels at $p$. In \cite{PR}, Pappas and Rapoport enhanced the moduli problem by adding the ``spin condition'' and conjectured its flatness. However, in \cite{smithling} Smithling gives a counterexample, and refines the moduli problem again with the ``strenghtened spin condition''. It has recently been proved in \cite{yu} (for a special case) and in \cite{luo} (for the general case) that Smithling's strenghtened spin model is flat. We note that the integral models of \cite{pappas, PR, smithling} are defined for arbitrary signature $(r,s)$ where $r+s=n$, however the flatness results mentioned above are dealing with the case of signature $(n-1,1)$ and $(1,n-1)$. Also, in loc. cit. the authors are working with local models of the Shimura varieties, but by the theory of linear modifications of \cite{pappas}, they also correspond to integral models of the actual Shimura varieties.\\
Let us go back to Pappas' integral model of \cite{pappas} with level at $p$ given by the stabilizer of a self-dual lattice and signature $(n-1,1)$. In loc. cit., it is proved that the integral model has isolated singularities located in the special fiber. Moreover, it is semi-stable when $n=2$, and its blow-up at the singular points is semi-stable when $n\geq 3$. The proof actually consists in proving the corresponding properties on the associated local model, which is a moduli space defined purely in terms of linear algebra, and that is represented by a closed subscheme of a grassmannian. The local model is related to the integral model of the Shimura variety via the local model diagram, which shows that both schemes are étale locally isomorphic. When trying to prove a property on a Shimura variety which is étale local, it is thus a common technique to reduce the problem to proving the corresponding statement on the local model, which is admittedly easier to work with. The local model constructed by Pappas has a single singular point located in its special fiber. By blowing-up at this singular point, one recovers Krämer's splitting model which was originally defined in \cite{kramer} by an enhanced moduli problem (see \cite{shi} Appendix A for a proof that the blow-up and the splitting model are isomorphic). Krämer showed that the special fiber of the splitting model consists of only two smooth divisors, one isomorphic to a projective space and the other having the structure of a $\mathbb P^1$-bundle over the scheme theoretic intersection, which itself is an explicit smooth projective quadric. We mention that splitting models for arbitrary maximal level at $p$ (ie. stabilizer of an almost self-dual lattice) have been introduced in \cite{HLS} and \cite{ZZ}, and that semi-stable reduction also holds in this context.\\
In this paper, we compute the cohomology sheaves $\mathrm R^i \Psi_{\eta} \Lambda$ of the $\ell$-adic complex of nearby cycles on Pappas' local model. To be more precise, let $p$ be a fixed odd prime which ramifies in $\mathbb E$, and let $E := \mathbb E_p$ denote the $p$-adic completion. Thus $E$ is a quadratic ramified extension of $\mathbb Q_p$ (in the body of the paper, $\mathbb Q_p$ is actually replaced with an arbitrary $p$-adic field $E_0$). Let $\pi$ denote a uniformizer of $E$ such that $\overline{\pi} = -\pi$, where $\overline{\,\cdot\,}$ denotes the non-trivial element of $\mathrm{Gal}(E/\mathbb Q_p)$. Let $M^{\mathrm{loc}}$ denote the local model associated to Pappas' flat integral model of the Shimura variety. Then $M^{\mathrm{loc}}$ is a flat projective scheme over $\mathrm{Spec}(\mathcal O_E)$, whose special fiber we denote by $\overline{M^{\mathrm{loc}}}$. The special fiber has dimension $n-1$ and it has a unique singular point $y^{\mathrm{sg}} \in \overline{M^{\mathrm{loc}}}(\mathbb F_p)$. Let $b:M^{\mathrm{bl}}\to M^{\mathrm{loc}}$ denote the blow-up of $M^{\mathrm{loc}}$ at $y^{\mathrm{sg}}$. Eventually, fix a prime $\ell \not = p$ and let $\Lambda$ denote a ring of coefficients (eg. $\Lambda = \mathbb Z / \ell^k \mathbb Z, \mathbb Z_{\ell}, \mathbb Q_{\ell}$ ; see Notations \ref{Notations} for the precise hypotheses). Let $\mathrm R\Psi_{\eta}\Lambda$ (resp. $\mathrm R\Psi_{\eta}^{\mathrm{bl}} \Lambda$) denote the associated nearby cycle complex on $\overline{M^{\mathrm{loc}}}$ (resp. on the blow-up $\overline{M^{\mathrm{bl}}}$), see \cite{illusie}. We will prove the following statement. 

\begin{theorem*}
We have 
\begin{equation*}
\mathrm R^i\Psi_{\eta}\Lambda = 
\begin{cases}
\Lambda & \text{if } i = 0,\\
i_{\overline{y^{\mathrm{sg}}}*}\Lambda[\epsilon p^{\frac{n}{2}}] & \text{if } n \text{ is even and } i = n-1,\\
0 & \text{else},
\end{cases}
\end{equation*}
where, for $n$ even, $i_{\overline{y^{\mathrm{sg}}}*}\Lambda[\epsilon p^{\frac{n}{2}}]$ is the skyscraper sheaf on $\overline{M^{\mathrm{loc}}}$ concentrated at the geometric singular point $\overline{y^{\mathrm{sg}}}$, with $\mathrm{Gal}(\overline E/E)$-action which is trivial on the inertia subgroup and such that the Frobenius element acts by multiplication by $\epsilon p^{\frac{n}{2}}$, where $\epsilon = \pm 1$.
\end{theorem*}

The sign $\epsilon$ is explicitly determined depending on the nature of the algebraic datum underlying the definition of $M^{\mathrm{loc}}$. Namely, it depends on whether a certain $E/\mathbb Q_p$-hermitian space is split or not. It is interesting to observe that, as a consequence of the computation above, we have $\mathrm R\Psi_{\eta}\Lambda \simeq \Lambda$ when $n$ is odd. Thus the nearby cycles are trivial despite the presence of a singularity. This fact was already proved by Zhu in \cite{zhu} Corollary 7.6, using a different and more general method.\\
The principle of the proof is simple, as it relies only on the proper base change formula applied to the blow-up morphism $b:M^{\mathrm{bl}} \to M^{\mathrm{loc}}$. Since $b$ induces an isomorphism on the generic fibers, we have 
\begin{equation*}
\mathrm R\Psi_{\eta}\Lambda \simeq \mathrm Rb_{\overline s *} \mathrm R\Psi_{\eta}^{\mathrm{bl}}\Lambda,
\end{equation*}
where $b_{\overline s}:(\overline{M^{\mathrm{bl}}})_{\overline{\mathbb F_p}} \to (\overline{M^{\mathrm{loc}}})_{\overline{\mathbb F_p}}$ denotes the morphism on geometric special fibers. Thus, if $\overline{x}$ is a point of $\overline{M^{\mathrm{loc}}}$, we have 
\begin{equation*}
(\mathrm R^i\Psi_{\eta}\Lambda)_{\overline x} \simeq \mathrm H^i(b_{\overline s}^{-1}\{\overline x\},\mathrm R\Psi_{\eta}^{\mathrm{bl}} \Lambda).
\end{equation*}
Clearly if $\overline x$ does not lie over $y^{\mathrm{sg}}$ then the stalk above vanishes for all $k\geq 1$. Denoting by $Z_1 := b_{s}^{-1}\{y^{\mathrm{sg}}\} \simeq \mathbb P^{n-1}$ the exceptional divisor, we are thus reduced to computing the cohomology of $Z_1$ with coefficients in the nearby cycles complex $\mathrm R\Psi_{\eta}^{\mathrm{bl}} \Lambda$. Since $M^{\mathrm{bl}}$ has semi-stable reduction, the cohomology sheaves of $\mathrm R\Psi_{\eta}^{\mathrm{bl}} \Lambda$ are well-known, cf. \cite{illusie}. We carry out the computation via diagram chasing through various Gysin sequences related to $\mathrm R\Psi_{\eta}^{\mathrm{bl}} \Lambda$ and its restriction to explicit subvarieties of $\overline{M^{\mathrm{bl}}}$. Well-known properties of nearby cycles, such as compatibility with duality and perversity, also play an important role.\\
We remark that Krämer has already computed the alternating semisimple trace of the Frobenius on the nearby cycles in \cite{kramer}. Namely, she computed the following sum
\begin{equation*}
\mathrm{Tr}^{\mathrm{ss}}(\mathrm{Frob}, \mathrm R\Psi_t(\mathbb Q_{\ell})_{y^{\mathrm{sg}}}) := \sum_i (-1)^i \mathrm{Tr}^{\mathrm{ss}}(\mathrm{Frob}, \mathrm R^i\Psi_t(\mathbb Q_{\ell})_{y^{\mathrm{sg}}}).
\end{equation*}
Here $\mathrm R^i\Psi_t \mathbb Q_{\ell} := (\mathrm R\Psi_{\eta} \mathbb Q_{\ell})^{P}$ denotes the tame nearby cycles, defined by taking the invariants under the largest $p$-subgroup $P\subset I$, and $I \subset \mathrm{Gal}(\overline E/E)$ denotes the inertia. Our computations recover Krämer's results and is more precise. Indeed, the inertia acting unipotently (even trivially), the semisimple trace agrees with the usual trace of the Frobenius. Moreover, the sum actually consists of only one or two non-zero terms depending on whether $n$ is odd or even. See Remark \ref{LastRemark} for more details.\\
As an application, our Theorem can be used to determine the cohomology sheaves of the nearby cycles on the actual Shimura variety. Let $\mathrm S$ denote an integral model over $\mathrm{Spec}(\mathcal O_E)$ of a Shimura variety whose local model is isomorphic to $M^{\mathrm{loc}}$. Typically, $\mathrm S$ can be defined via the moduli problem of \cite{pappas} Definition 3.4. Then $\mathrm S$ and $M^{\mathrm{loc}}$ are related via a local model diagram in the sense of \cite{pappas} Theorem 2.2. Let $\mathrm R\Psi_{\eta}^{\mathrm S}\Lambda$ denote the nearby cycles on the Shimura variety. Just as $M^{\mathrm{loc}}$, the Shimura variety $\mathrm S$ has isolated singularities located in its special fiber $\overline{\mathrm S}$. It readily follows that $\mathrm R^0\Psi_{\eta}^{\mathrm S}\Lambda = \Lambda$ and that the higher cohomology sheaves $\mathrm R^i\Psi_{\eta}^{\mathrm S}\Lambda$ for $i>0$ are skyscraper, concentrated on the singular points. To any point $\overline x \in \overline{\mathrm S}(\overline{\mathbb F_p})$, there corresponds a point $\overline y \in \overline{M^{\mathrm{loc}}}(\overline{\mathbb F_p})$. Moreover, $\overline{x}$ is singular if and only if $\overline y = \overline{y^{\mathrm{sg}}}$. Through the local model diagram, the stalk $(\mathrm R^i\Psi_{\eta}^{\mathrm S}\Lambda)_{\overline x}$ is isomorphic to $(\mathrm R^i\Psi_{\eta}\Lambda)_{\overline y}$, which is computed by our Theorem. 

\subsection*{Acknowledgements}

The author is grateful to the referees for their careful reading of the manuscript, for numerous helpful suggestions, and for pointing out a mistake in a key argument of an earlier version of this paper. 

\section{Notations} \label{Notations}

In this paper, we fix a prime number $\ell$ and a coefficient ring $\Lambda$ which is one of the following kinds
\begin{itemize}
\item $\Lambda = K$ where $K$ is any algebraic extension of $\mathbb Q_{\ell}$,
\item $\Lambda = \mathcal O_K$ where $K$ is any algebraic extension of $\mathbb Q_{\ell}$,
\item $\Lambda$ is a regular ring killed by a power of $\ell$ (here, we mean ``regular'' in the sense of \citep{HansenScholze}, ie. any truncation of a perfect complex of $\Lambda$-modules is still perfect).
\end{itemize}
Given $f: X \to \mathrm{Spec}(k)$ a variety (integral, separated and of finite type) over an algebraically closed field with $\mathrm{char}(k) \not = \ell$, we write $D_c^b(X,\Lambda)$ for the bounded ``derived'' category of constructible $\Lambda$-sheaves on $X$. For $K \in D_c^b(X,\Lambda)$, we write $D_X(K) := \mathrm R \mathcal Hom(K,f^!\Lambda)$ for its Verdier dual. It defines a functor $D_X: D_c^b(X,\Lambda) \to D_c^b(X,\Lambda)$. When the context is clear, we write $D$ instead of $D_X$. We denote by $\mathrm{Perv}(X,\Lambda) \subset D_c^b(X,\Lambda)$ the abelian category of perverse sheaves on $X$, cf. \cite{HansenScholze}.\\
Whenever we work with a variety $X$ over a field $k$ of characteristic different from $\ell$, not necessarily algebraically closed, we always write $\mathrm H^{\bullet}(X,\Lambda)$ instead of $\mathrm H^{\bullet}(X \otimes_k \overline{k},\Lambda)$. 

\section{Geometry of the local model}\label{Section1}

Let $p$ be an odd prime, let $E_0$ be a $p$-adic field with ring of integers $\mathcal O_{E_0}$, uniformizer $\pi_0$ and residue field $\kappa = \mathcal O_{E_0}/(\pi_0) \simeq \mathbb F_q$ for some $q$ which is a power of $p$. Let $E$ be a ramified quadratic extension of $E_0$ with ring of integers $\mathcal O_E$. Fix a uniformizer $\pi$ of $E$ such that $\overline{\pi} = -\pi$ where $\overline{\,\cdot\,} \in \mathrm{Gal}(E/E_0)$ is the non-trivial Galois involution on $E$. Let $V$ be an $E$-vector space of dimension $n\geq 2$ equipped with a perfect $E/E_0$-hermitian form $(\cdot,\cdot):V\times V \to E$. By convention, we assume that it is linear in the first variable and semilinear in the second. We also define an $E_0$-valued alternating bilinear form $\psi$ via the formula 
\begin{equation*}
\psi(u,v) := \mathrm{Tr}_{E/E_0}\left(\frac{(u,v)}{2\pi}\right),
\end{equation*}
for all $u,v\in V$. Note that one may also recover $(\cdot,\cdot)$ from $\psi$. Let $L$ be a self-dual (for $(\cdot,\cdot)$ or equivalently for $\psi$) $\mathcal O_E$-lattice in $V$. The local model of signature $(n-1,1)$ considered by Pappas in \cite{pappas} is defined as follows. For an $\mathcal O_{E}$-scheme $S$, let $M^{\psi}_{n-1,1}(S)$ denote the set of subsheaves $\mathcal F \subset L \otimes_{\mathcal O_{E_0}} \mathcal O_S$ of $\mathcal O_E \otimes_{\mathcal O_{E_0}} \mathcal O_S$-modules such that 
\begin{enumerate}[label={--},noitemsep]
\item as a sheaf of $\mathcal O_S$-modules, $\mathcal F$ is locally a free direct summand of $L \otimes_{\mathcal O_{E_0}} \mathcal O_S$ of rank $n$,
\item $\mathcal F$ is totally isotropic for the alternating bilinear form $\psi \otimes \mathcal O_S$ on $L \otimes_{\mathcal O_{E_0}} \mathcal O_S$,
\item the Kottwitz and Pappas conditions are satisfied:
\begin{align*}
\forall x \in \mathcal O_{E},\, & \det(T-x\,|\,\mathcal F) = (T-x)^{1}(T-\overline{x})^{n-1} \in \mathcal O_{E}[T],\\
\text{if } n\geq 3, & \bigwedge^{n}(x\otimes 1 - 1\otimes x) = 0 \text{ and } \bigwedge^{2}(x\otimes 1 - 1\otimes \overline{x}) = 0 \text{ on } \mathcal F.
\end{align*}
\end{enumerate}
This functor is represented by a projective $\mathcal O_E$-scheme.

\begin{remark}
In \cite{kramer}, the moduli space above is denoted $M_{n-1,1}^{\mathrm{loc}}$, and in \cite{pappas} it is denoted $M^{\prime\psi}_{n-1,1}$. The definition in \cite{pappas} is slightly different, in so far as the Kottwitz and Pappas conditions are imposed on the quotient $\mathcal Q := L\otimes_{\mathcal O_{E_0}} \mathcal O_S / \mathcal F$ instead of $\mathcal F$. However, both definitions coincide as it is not difficult to check that these conditions on $\mathcal Q$ are equivalent to the same conditions on $\mathcal F$ with the roles of $n-1$ and $1$ reversed. Also, in \cite{pappas} only the case $E_0 = \mathbb Q_p$ is considered, however in \cite{kramer} $E_0$ is allowed to be any $p$-adic field.
\end{remark}

From now on, we write $M^{\psi} := M_{n-1,1}^{\psi}$ as a scheme over $\mathrm{Spec}(\mathcal O_E)$ (corresponding to $M^{\mathrm{loc}}$ in the introduction). According to \cite{pappas}, $M^{\psi}$ is flat over $\mathrm{Spec}(\mathcal O_E)$, regular if $n = 2$, and regular outside of a single closed point $y^{\mathrm{sg}}$ if $n\geq 3$. Explicitly, $y^{\mathrm{sg}}$ lies in the special fiber of $M^{\psi}$. It is the $\mathbb F_q$-rational point given by $\mathcal F^{\mathrm{sg}} := (\pi\otimes 1) L \otimes_{\mathcal O_{E_0}} \mathbb F_q$. If $n\geq 3$, let $b:M^{\psi,\mathrm{bl}}\to M^{\psi}$ be the blow-up of $M^{\psi}$ along $y^{\mathrm{sg}}$. Let $\overline{M^{\psi}}$ and $\overline{M^{\psi,\mathrm{bl}}}$ denote the special fibers over $\kappa = \mathbb F_q$ of $M^{\psi}$ and of $M^{\psi,\mathrm{bl}}$ respectively. If $n=2$, $\overline{M^{\psi}}$ is a divisor with simple normal crossings in $M^{\psi}$. If $n\geq 3$, then $M^{\psi,\mathrm{bl}}$ is regular and $\overline{M^{\psi,\mathrm{bl}}}$ is a divisor with simple normal crossings in $M^{\psi,\mathrm{bl}}$. In \cite{kramer}, the author gives an explicit description of the geometry of the special fiber $\overline{M^{\psi,\mathrm{bl}}}$. 

\begin{theorem}
The special fiber $\overline{M^{\psi,\mathrm{bl}}}$ is the union of two smooth irreducible varieties $Z_1$ and $Z_2$. We have $Z_1 := b^{-1}\{y^{\mathrm{sg}}\} \simeq \mathbb P^{n-1}$ and $Z_2$ is a $\mathbb P^1$-bundle over the scheme theoretic intersection $Q := Z_1\cap Z_2$. Moreover, the closed immersion $\iota_1:Q\hookrightarrow Z_1$ identifies $Q$ with an explicit smooth quadric in $\mathbb P^{n-1}$.
\end{theorem}

We may sum up the situation via the following diagram

\begin{center}
\begin{tikzcd}[ampersand replacement=\&]
	\& \overline{M^{\psi,\mathrm{bl}}} \& \\
	Z_1 \simeq \mathbb P^{n-1} \arrow[hook]{ur}{i_1} \& \& Z_2 \simeq \mathbb P(\mathcal E) \arrow[hook',swap]{ul}{i_2} \arrow[swap, bend right=30]{dl}{\mathrm{pr}} \\
	\& Q \arrow[hook']{ul}{\iota_1} \arrow[hook,swap]{ur}{\iota_2} \&
\end{tikzcd}
\end{center}

\noindent where $\mathcal E$ is a certain locally free sheaf of rank $2$ on $Q$. We write $i_Q := i_1\circ \iota_1 = i_2 \circ \iota_2$ for the closed immersion $Q \hookrightarrow \overline{M^{\psi,\mathrm{bl}}}$. Moreover, we also define $U_i := \overline{M^{\psi,\mathrm{bl}}} \setminus Z_i = Z_{i+1} \setminus Q$, where the index $i$ is seen in $\mathbb Z/2\mathbb Z$. We write $j_i:U_i \hookrightarrow \overline{M^{\psi,\mathrm{bl}}}$ and $j_i': U_i \hookrightarrow Z_{i+1}$ for the open immersions.\\
Let us give more details. The blow-up $M^{\psi,\mathrm{bl}}$ can also be characterized as the moduli space classifying pairs $(\mathcal F_0,\mathcal F)$ of $\mathcal O_E \otimes_{\mathcal O_{E_0}} \mathcal O_S$-submodules of $L \otimes_{\mathcal O_{E_0}} \mathcal O_S$ (for $S$ an $\mathcal O_E$-scheme) such that 
\begin{enumerate}[label={--},noitemsep]
\item $\mathcal F \in M^{\psi}(S)$,
\item as a sheaf of $\mathcal O_S$-modules, $\mathcal F_0$ is a locally free direct summand of $L \otimes_{\mathcal O_{E_0}} \mathcal O_S$ of rank $1$,
\item $\mathcal F_0 \subset \mathcal F$,
\item $(\pi\otimes 1 + 1\otimes \pi)\mathcal F \subset \mathcal F_0$,
\item $(\pi\otimes 1 - 1\otimes \pi) \mathcal F_0 = 0.$
\end{enumerate}

The blow-up morphism $b:M^{\psi,\mathrm{bl}} \to M^{\psi}$ corresponds to the forgetful functor $(\mathcal F,\mathcal F_0) \mapsto \mathcal F$.

\begin{remark}
Both papers \cite{pappas} and \cite{kramer} study $M^{\psi,\mathrm{bl}}$, but the former uses the blow-up construction and the latter used the moduli description. The equivalence between both definitions was well-known for a long time, and a proof may be found in \cite{shi} Appendix A.
\end{remark}

The map $\mathrm{pr}: \overline{M^{\psi,\mathrm{bl}}} \to \mathbb P((\pi\otimes 1)L\otimes_{\mathcal O_{E_0}}\mathbb F_q) \simeq \mathbb P^{n-1}$ induced by $(\mathcal F,\mathcal F_0) \mapsto \mathcal F_0$ restricts to an isomorphism on the fiber $Z_1 := b^{-1}\{y^{\mathrm{sg}}\}$. Krämer then defines a pairing on $(\pi\otimes 1)L\otimes_{\mathcal O_{E_0}}\mathbb F_q$ by the formula 
\begin{equation*}
\{ (\pi\otimes 1)v,(\pi\otimes 1)w\} := \psi((\pi\otimes 1)v,w).
\end{equation*}
It is easy to check that it is well-defined, symmetric and non-degenerate. The closed subvariety $Q \subset Z_1$ consists of the isotropic lines in $(\pi\otimes 1)L\otimes_{\mathcal O_{E_0}}\mathbb F_q$. 

\begin{proposition}
The variety $Q$ is smooth and geometrically irreducible of dimension $n-2$.
\end{proposition}

\begin{proof}
Let $\delta \in \mathcal O_{E_0}^{\times}$ which is not the norm of any element of $\mathcal O_E$. There is an $E$-basis $(e_1,\ldots ,e_n)$ of $V$ in which $(\cdot,\cdot)$ is given by the matrix $\mathrm{Diag}(1,\ldots,1)$ or $\mathrm{Diag}(1,\ldots,1,\delta)$. We may further assume that $L$ is generated by the $e_i$'s. A line of $(\pi\otimes 1)L\otimes_{\mathcal O_{E_0}} \mathbb F_q$ generated by a non-zero vector $(\pi\otimes 1)(a_1e_1 + \ldots + a_ne_n)$ is then isotropic if and only if $\sum_{i=1}^{n} a_i^2 = 0$ or $\sum_{i=1}^{n-1} a_i^2 + \overline{\delta} a_n^2 = 0$, where $\overline{\delta} \in \mathbb F_q$ is the residue modulo $\pi_0$ of $\delta$. Since $\mathbb F_q$ has odd characteristic, the proposition follows.
\end{proof}

Let $Z_2 := \mathrm{pr}^{-1}(Q)$. Then Krämer builds a locally free sheaf $\mathcal E$ of rank $2$ on $Q$ and an isomorphism $Z_2 \simeq \mathbb P(\mathcal E)$. Eventually, one checks that $Q$ coincides with the scheme theoretical intersection $Z_1\cap Z_2$.

\section{Auxiliary cohomological computations} \label{Section2}

We explicitly determine the cohomology of all the subvarieties of $\overline{M^{\psi,\mathrm{bl}}}$ introduced in Section \ref{Section1}, as it will be useful towards the computation of the nearby cycles.

\begin{proposition}\label{CohomologyOfQ}
For $0 \leq i \leq 2(n-2)$, the restriction map $\iota_1^*: \mathrm H^i(Z_1,\Lambda) \to \mathrm H^i(Q,\Lambda)$ is an isomorphism if $i \not = n-2$, and is injective if $i = n-2$. The middle degree cohomology group decomposes naturally as 
\begin{equation*}
\mathrm H^{n-2}(Q,\Lambda) \simeq \iota_1^*(\mathrm H^{n-2}(Z_1,\Lambda))\oplus \mathrm H^{n-2}_{\mathrm{prim}}(Q,\Lambda),
\end{equation*}
and we have
\begin{equation*}
\mathrm H^{n-2}_{\mathrm{prim}}(Q,\Lambda) = \mathrm{Ker}(\iota_{1*}:\mathrm H^{n-2}(Q,\Lambda) \to \mathrm H^n(Z_1,\Lambda)(1)) =
\begin{cases}
0 & \text{if } n \text{ is odd},\\
\Lambda[\epsilon q^{\frac{n-2}{2}}] & \text{if } n \text{ is even},
\end{cases}
\end{equation*}
where $\Lambda[\epsilon q^{\frac{n-2}{2}}]$ is the module $\Lambda$ with Frobenius action given by multiplication by $\epsilon q^{\frac{n-2}{2}}$, and $\epsilon = 1$ (resp. $\epsilon = -1$) if $(V,(\cdot,\cdot))$ is split (resp. non-split).
\end{proposition}

In particular, we have $\mathrm H^i(Q,\Lambda) = 0$ whenever $i$ is odd.

\begin{proof}
The first statement on the restriction map $\iota_1^*$ and the decomposition 
\begin{equation*}
\mathrm H^{n-2}(Q,\Lambda) \simeq \iota_1^*(\mathrm H^{n-2}(Z_1,\Lambda))\oplus \mathrm H^{n-2}_{\mathrm{prim}}(Q,\Lambda),
\end{equation*}
where $\mathrm H^{n-2}_{\mathrm{prim}}(Q,\Lambda) := \mathrm{Ker}(\iota_{1*}:\mathrm H^{n-2}(Q,\Lambda) \to \mathrm H^n(Z_1,\Lambda)(1))$ is the kernel of the Gysin map, follow from the weak Lefschetz theorem and Poincaré duality. Moreover, the primitive cohomology $\mathrm H^{n-2}_{\mathrm{prim}}(Q,\Lambda)$ is a free $\Lambda$-module whose rank depends only on the degree of the hypersurface $Q$, see \cite{milne} Example 16.4. Since $Q$ is a smooth quadric in $Z_1 \simeq \mathbb P^{n-1}$, we know that $\mathrm H^{n-2}_{\mathrm{prim}}(Q,\Lambda) = 0$ if $n$ is odd and that $\mathrm{rank} (\mathrm H^{n-2}_{\mathrm{prim}}(Q,\Lambda)) = 1$ if $n$ is even. Thus, the non-trivial part of the Proposition consists in the computation of the Frobenius eigenvalue on the primitive cohomology when $n = 2m \geq 4$ is even, which we now assume until the end of the proof. By the Lefschetz trace formula, see \cite{sga41/2} ``Rapport sur la formule des traces'', we have 
\begin{equation*}
\# Q(\mathbb F_q)\cdot 1_{\Lambda} = \sum_{i\geq 0} (-1)^i\mathrm{Trace}(\mathrm{Frob}\,|\, \mathrm H^{i}(Q,\Lambda)) = 1 + q + \ldots + q^{n-2} + \mathrm{Trace}(\mathrm{Frob} \,|\, \mathrm H^{n-2}_{\mathrm{prim}}(Q,\Lambda)).
\end{equation*}

\begin{remark} 
The symbol $1_{\Lambda}$ is the unity element of $\Lambda$. The formalism developped in \cite{HansenScholze} guarantees that the quantities involed in the formula make sense under the hypotheses on $\Lambda$ imposed in \ref{Notations}.
\end{remark}

The left-hand side can be computed using \cite{weil}. To do this, we need to fix an equation for $Q$. The hermitian space $V$ is split if and only if $(\cdot,\cdot)$ is given by the matrix $\mathrm{Diag}(1,\ldots, 1, -1, \ldots, -1)$ in some basis, where $1$ and $-1$ occur $m$ times each. It is non-split if and only if $(\cdot,\cdot)$ is given by $\mathrm{Diag}(1,\ldots ,1, -1,\ldots ,-1,-\delta)$ for some $\delta \in \mathcal O_{E_0}^{\times}$ which is not the norm of any unit of $\mathcal O_E$. Assuming that $L$ is generated by such a basis, the quadric $Q \subset \mathbb P^{n-1}$ is given by the equation $x_1^2 + \ldots + x_m^2 - x_{m+1}^2 - \ldots - x_n^2 = 0$ in the split case, and by $x_1^2 + \ldots + x_m^2 - x_{m+1}^2 - \ldots - x_{n-1}^2 - \overline{\delta} x_n^2 = 0$ in the non-split case, where $\overline{\delta} \in \mathbb F_q$ is the residue modulo $\pi$ of $\delta$. Consider the following Jacobi sum 
\begin{equation*}
j_m := \frac{1}{q-1} \sum_{\substack{u_1 + \ldots + u_{2m} = 0\\ u_i \in \mathbb F_q^{\times}}} \left(\frac{u_1}{q}\right) \ldots \left(\frac{u_{2m}}{q}\right),
\end{equation*}
where $\left(\frac{\,\cdot\,}{q}\right)$ denotes the Legendre symbol on $\mathbb F_q$. According to \cite{weil}, we have 
\begin{equation*}
\#Q(\mathbb F_q) = \begin{cases}
1 + \ldots + q^{n-2} + \left(\frac{-1}{q}\right)^m j_m & \text{if } (V,(\cdot,\cdot)) \text{ is split},\\
1 + \ldots + q^{n-2} + \left(\frac{-1}{q}\right)^m\left(\frac{\overline{\delta}}{q}\right)^{-1} j_m & \text{if } (V,(\cdot,\cdot)) \text{ is non-split}.
\end{cases}
\end{equation*}
Since $\overline{\delta}$ is not a square in $\mathbb F_q$ (otherwise, since $p\not = 2$, $\delta$ would be the square of an element $\alpha \in \mathcal O_{E_0}$ by Hensel's lemma, so that we would have $\mathrm{Norm}_{E/E_0}(\alpha) = \alpha^2 = \delta$, a contradiction), we have $\left(\frac{\overline{\delta}}{q}\right) = -1$. Thus, the proof is over once we compute $j_m$, which is the object of the next Lemma.
\end{proof}

\begin{lemma}
We have
\begin{equation*}
j_m = \left(\frac{-1}{q}\right)^m q^{m-1}.
\end{equation*}
\end{lemma}

\begin{proof}
We establish a relation between $j_m$ and $j_{m-1}$. Let us assume that $m\geq 2$. Observe that the equation $u_1 + \ldots + u_{2m} = 0$ is equivalent to $v_1 + \ldots + v_{2m-1} = 1$ where $v_i := -\frac{u_i}{u_{2m}}$. Thus, we rearrange the sum as 
\begin{equation*}
(q-1)j_m = \sum_{u_{2m} \in \mathbb F_q^{\times}} \left(\frac{-1}{q}\right)^{2m-1}\left(\frac{u_{2m}}{q}\right)^{2m} \sum_{\substack{v_1 + \ldots + v_{2m-1} = 1\\ v_i \in \mathbb F_q^{\times}}} \left(\frac{v_1}{q}\right) \ldots \left(\frac{v_{2m-1}}{q}\right). 
\end{equation*}
Since the Legendre symbol takes value in $\{\pm 1\}$, this simplifies to
\begin{equation}\label{eq1}
j_m = \left(\frac{-1}{q}\right)\sum_{\substack{v_1 + \ldots + v_{2m-1} = 1\\ v_i \in \mathbb F_q^{\times}}} \left(\frac{v_1}{q}\right) \ldots \left(\frac{v_{2m-1}}{q}\right).
\end{equation}
Note that the factor $(q-1)$ on the LHS simplifies with the sum over $u_{2m} \in \mathbb F_{q}^{\times}$, as the summand of the sum does not actually depend on $u_{2m}$ anymore. By the change of variable $w := 1 - v_{2m-1} \in \mathbb F_q \setminus\{1\}$, we have
\begin{equation*}
j_m = \left(\frac{-1}{q}\right)\sum_{\substack{v_1 + \ldots + v_{2m-2} = w\\ v_i \in \mathbb F_q^{\times} \\ w \in \mathbb F_q \setminus\{1\}}} \left(\frac{v_1}{q}\right) \ldots \left(\frac{v_{2m-2}}{q}\right) \left(\frac{1-w}{q}\right).\end{equation*}
Isolating the terms corresponding to $w=0$, we obtain
\begin{equation}\label{eq2}
j_m = \left(\frac{-1}{q}\right)(q-1)j_{m-1} + \left(\frac{-1}{q}\right)\sum_{\substack{v_1 + \ldots + v_{2m-2} = w\\ v_i \in \mathbb F_q^{\times} \\ w \in \mathbb F_q \setminus\{0,1\}}} \left(\frac{v_1}{q}\right) \ldots \left(\frac{v_{2m-2}}{q}\right) \left(\frac{1-w}{q}\right).
\end{equation}
Let us call $S$ the sum on the right. For $1\leq i \leq 2m-2$ we write $w_i := \frac{v_i}{w}$, so that 
\begin{equation*}
S = \sum_{w \in \mathbb F_q \setminus \{0,1\}} \left(\frac{1-w}{q}\right)\left(\frac{w}{q}\right)^{2m-2}\sum_{\substack{w_1+\ldots +w_{2m-2} = 1\\ w_i \in \mathbb F_q^{\times}}} \left(\frac{w_1}{q}\right) \ldots \left(\frac{w_{2m-2}}{q}\right).
\end{equation*}
Since $2m-2$ is even and since $\sum_{w\in \mathbb F_q \setminus \{0,1\}} \left(\frac{1-w}{q}\right) = -1$, this simplifies to 
\begin{equation*}
S = -\sum_{\substack{w_1+\ldots +w_{2m-2} = 1\\ w_i \in \mathbb F_q^{\times}}} \left(\frac{w_1}{q}\right) \ldots \left(\frac{w_{2m-2}}{q}\right).
\end{equation*}
If $m = 2$, then $S = -\sum_{x\in \mathbb F_q \setminus\{0,1\}} \left(\frac{x(1-x)}{q}\right)$. Notice that $\left(\frac{x(1-x)}{q}\right) = \left(\frac{x^{-1}-1}{q}\right)$. Since $x\mapsto x^{-1}$ defines a bijection on $\mathbb F_q \setminus \{0,1\}$, by a change of variable we obtain 
\begin{equation*}
S = -\sum_{y \in \mathbb F_q\setminus\{0,1\}} \left(\frac{y-1}{q}\right) = \left(\frac{-1}{q}\right) = j_1.
\end{equation*}
For the last equality, we can indeed compute $j_1$ directly as follows
\begin{equation*}
j_1 = \frac{1}{q-1}\sum_{x \in \mathbb F_q^{\times}} \left(\frac{x}{q}\right)\left(\frac{-x}{q}\right) = \frac{1}{q-1}\sum_{x \in \mathbb F_q^{\times}} \left(\frac{-1}{q}\right) = \left(\frac{-1}{q}\right).
\end{equation*}
Now assume that $m \geq 3$. We introduce the change of variable $x := 1 - w_{2m-2}$ and write 
\begin{equation*}
S = -\sum_{\substack{w_1 + \ldots + w_{2m-3} = x \\ w_i \in \mathbb F_q^{\times} \\ x \in \mathbb F_q \setminus\{1\}}} \left(\frac{w_1}{q}\right) \ldots \left(\frac{w_{2m-3}}{q}\right)\left(\frac{1-x}{q}\right).
\end{equation*}
Isolating the terms corresponding to $x=0$, we have 
\begin{equation*}
S = -\sum_{\substack{w_1 + \ldots + w_{2m-3} = 0 \\ w_i \in \mathbb F_q^{\times}}} \left(\frac{w_1}{q}\right) \ldots \left(\frac{w_{2m-3}}{q}\right) - \sum_{\substack{w_1 + \ldots + w_{2m-3} = x \\ w_i \in \mathbb F_q^{\times} \\ x \in \mathbb F_q \setminus\{0,1\}}}\left(\frac{w_1}{q}\right) \ldots \left(\frac{w_{2m-3}}{q}\right)\left(\frac{1-x}{q}\right).
\end{equation*}
First we compute the sum on the left. Fix an element $\lambda \in \mathbb F_q^{\times}$ which is not a square, and make the change of variables $w_i \mapsto \lambda w_i$. Since $\left(\frac{\lambda}{q}\right)^{2m-3} = -1$, one finds out that this sum vanishes. Next, we proceed to the change of variables $x_i := \frac{w_i}{x}$ in the sum on the right. We obtain 
\begin{equation*}
S = -\sum_{x\in \mathbb F_q\setminus\{0,1\}} \left(\frac{1-x}{q}\right) \left(\frac{x}{q}\right)^{2m-3}\sum_{\substack{x_1 + \ldots + x_{2m-3} = 1 \\ x_i \in \mathbb F_q^{\times}}} \left(\frac{x_1}{q}\right) \ldots \left(\frac{x_{2m-3}}{q}\right).
\end{equation*}
Since $2m-3$ is odd, the sum on $x$ recovers the case $m=2$ treated above, and simplifies to $\left(\frac{-1}{q}\right)$. On the other hand, according to \eqref{eq1} with $m$ replaced by $m-1 \geq 2$, the inner sum is related to $j_{m-1}$. All in all, we obtain 
\begin{equation*}
S = j_{m-1}.
\end{equation*}
Replacing in \eqref{eq2}, we deduce that 
\begin{equation*}
j_m = \left(\frac{-1}{q}\right)(q-1)j_{m-1} + \left(\frac{-1}{q}\right)j_{m-1} = \left(\frac{-1}{q}\right)qj_{m-1},\end{equation*}
and the proof is over.
\end{proof}

\begin{lemma}\label{CohomologyZ2}
The pullback $\mathrm{pr}^{*}$ gives an isomorphism of graded $\mathrm H^{*}(Q,\Lambda)$-modules 
\begin{equation*}
\mathrm H^{*}(Q,\Lambda)[t]/(t^2) \xrightarrow{\sim} \mathrm H^{*}(Z_2,\Lambda),
\end{equation*}
sending $t$ to the image $\zeta := [\mathcal O_{\mathbb P(\mathcal E)}(1)] \in \mathrm H^2(Z_2,\Lambda)(1)$ of $\mathcal O_{\mathbb P(\mathcal E)}(1)$ via the Kummer sequence.
\end{lemma}

\begin{proof}
This is the projective bundle formula, see for instance \cite{milne} Theorem 23.2.
\end{proof}

\begin{lemma}\label{CohomologyU1}
The restriction of $\mathrm{pr}$ to $U_1$ induces an isomorphism $(\mathrm{pr}_{|U_1})^{*}:\mathrm H^{i}(Q,\Lambda) \xrightarrow{\sim} \mathrm H^{i}(U_1,\Lambda)$ for all $i$. 
\end{lemma}

\begin{proof}
This is clear since $\mathrm{pr}_{|U_1}:U_1 \to Q$ is an $\mathbb A^1$-bundle on $Q$. 
\end{proof}

\begin{lemma}\label{CohomologyU2}
The restriction $(j_2')^*$ induces an isomorphism $\mathrm H^0(Z_1,\Lambda) \xrightarrow{\sim} \mathrm H^{0}(U_2,\Lambda)$. When $n$ is even, the connecting morphism of the Gysin sequence induces an isomorphism $\mathrm H^{n-1}(U_2,\Lambda) \xrightarrow{\sim} \mathrm H^{n-2}_{\mathrm{prim}}(Q,\Lambda)(-1)$. For all other $i$, we have $\mathrm H^i(U_2,\Lambda) = 0$.
\end{lemma}

\begin{proof}
The Gysin sequence associated to $\iota_1: Q \hookrightarrow Z_1 \hookleftarrow U_2 : j_2'$ yields, for all $i \in \mathbb Z$, an exact sequence
\begin{equation*}
0 \to \mathrm H^{2i-1}(U_2,\Lambda) \to \mathrm H^{2i-2}(Q,\Lambda)(-1) \xrightarrow{\iota_{1*}} H^{2i}(Z_1,\Lambda) \xrightarrow{(j_2')^*} \mathrm H^{2i}(U_2,\Lambda) \to 0.
\end{equation*}
By Proposition \ref{CohomologyOfQ}, the Gysin map $\iota_{1*}$ is an isomorphism for all $1 \leq i \leq n-1$, except when $n$ is even and $i = \frac{n}{2}$. In this case, it is surjective with kernel equal to $\mathrm H^{n-2}_{\mathrm{prim}}(Q,\Lambda)(-1)$. The result follows.
\end{proof}

\section{The nearby cycles on the local model}

%Recall from \cite{illusie} Paragraph 3.8 or \cite{saito} Section 2.2 that $\mathrm R \Psi_{\eta}\Lambda$ is an object of the abelian category of $-(n-1)$-shifted perverse sheaves $\mathrm{Perv}((\overline{M^{\psi,\mathrm{bl}}})_{\overline{\mathbb F_q}}, \Lambda)[-(n-1)] \subset D_c^b((\overline{M^{\psi,\mathrm{bl}}})_{\overline{\mathbb F_q}}, \Lambda)$, where the right-hand side denotes the bounded ``derived'' category of constructible $\Lambda$-sheaves. Inside $\mathrm{Perv}(\overline{M^{\psi,\mathrm{bl}}})_{\overline{\mathbb F_q}}, \Lambda)[-(n-1)]$, the nearby cycles admit a monodromy filtration 
%$$\ldots \subset F^i\mathrm R \Psi_{\eta}^{\mathrm{bl}}\Lambda \subset F^{i+1}\mathrm R \Psi_{\eta}^{\mathrm{bl}}\Lambda \subset \ldots,$$
%such that the graded pieces are given by (\cite{saito} Corollary 2.8)
%$$\mathrm{Gr}_{r}\mathrm R \Psi_{\eta}^{\mathrm{bl}}\Lambda \simeq \begin{cases}
%i_{Q*}\Lambda[-1] & \text{if } r = -1,\\
%a_{0*}\Lambda[0] & \text{if } r = 0,\\
%i_{Q*}\Lambda(-1)[-1] & \text{if } r = 1,\\
%0 & \text{else},
%\end{cases}$$
%where $a_0 : Z_1 \sqcup Z_2 \to \overline{M^{\psi,\mathrm{bl}}}$ is the natural map. It is easy to check that 
%\begin{equation*}
%F^i\mathrm R \Psi_{\eta}^{\mathrm{bl}}\Lambda = \begin{cases}
%0 & \text{if } r < -1,\\
%i_{Q*}\Lambda[-1] & \text{if } r = -1,\\
%\Lambda[0] & \text{if } r = 0,\\
%\mathrm R \Psi_{\eta}^{\mathrm{bl}}\Lambda & \text{if } r \geq 1.
%\end{cases}
%\end{equation*}

Let $\ell$ be a prime different from $p$ and let $\Lambda$ be a coefficient ring as in the Notations \ref{Notations}. Let $\mathrm R\Psi_{\eta} \Lambda$ denote the nearby cycles on the local model $M^{\psi}$. The rest of the paper is dedicated to proving the following statement.

%By the local model diagram, if $\overline x$ is a geometric point of $\overline{\mathrm S}_{K^p}$ and $\overline y = qs(\overline x)$ is the corresponding point of $\overline M$, then we have an isomorphism $(\mathrm R\Psi_{\eta}\Lambda)_{\overline x} \simeq (\mathrm R\Psi_{\eta}^{M}\Lambda)_{\overline y}$. Therefore the nearby cycles can be computed on the local model. 

\begin{theorem}\label{ComputationOfNearbyCycles}
We have 
\begin{equation*}
\mathrm R^i\Psi_{\eta}\Lambda = 
\begin{cases}
\Lambda & \text{if } i = 0,\\
i_{\overline{y^{\mathrm{sg}}}*}\Lambda[\epsilon q^{\frac{n}{2}}] & \text{if } n \text{ is even and } i = n-1,\\
0 & \text{else},
\end{cases}
\end{equation*}
where, for $n$ even, $i_{\overline{y^{\mathrm{sg}}}*}\Lambda[\epsilon q^{\frac{n}{2}}]$ is the skyscraper sheaf concentrated at the geometric singular point $\overline{y^{\mathrm{sg}}}$ with $\mathrm{Gal}(\overline E/E)$-action which is trivial on the inertia subgroup, and with Frobenius action given by multiplication by $\epsilon q^{\frac{n}{2}}$, where $\epsilon = 1$ if $n=2$ or if $n \geq 4$ and the hermitian space $(V,(\cdot,\cdot))$ is split, and $\epsilon = -1$ if $n \geq 4$ and the hermitian space is non-split.
\end{theorem}

Recall that $(V,(\cdot,\cdot))$ is said to be split if its discriminant $\mathrm{disc}(V) := (-1)^{\frac{n(n-1)}{2}}\mathrm{det}(V) \in E_0^{\times} / \mathrm{Norm}_{E/E_0}(E^{\times})$ is trivial, and non-split otherwise. The case $n=2$ is easy since $M^{\psi}$ already has semi-stable reduction. Thus in the remaining of this paper, we assume that $n\geq 3$. Let $\mathrm R \Psi_{\eta}^{\mathrm{bl}}\Lambda$ denote the nearby cycles on the blow-up $M^{\psi,\mathrm{bl}}$ of the local model. Since $M^{\psi,\mathrm{bl}}$ has semi-stable reduction and the special fiber $\overline{M^{\psi,\mathrm{bl}}}$ has only two irreducible components, we have 
\begin{equation*}
\mathrm R^i \Psi_{\eta}^{\mathrm{bl}}\Lambda = \begin{cases}
\Lambda & \text{if } i=0,\\
i_{Q*}\Lambda(-1) & \text{if } i=1,\\
0 & \text{else},
\end{cases}
\end{equation*}
see \cite{illusie} Théorème 3.2. By proper base change and since $b$ is an isomorphism on the generic fibers, we have
\begin{equation*}
\mathrm R \Psi_{\eta}\Lambda \simeq \mathrm R b_{\overline s*}\mathrm R \Psi_{\eta}^{\mathrm{bl}}\Lambda,
\end{equation*}
where $b_{\overline s}:(\overline{M^{\psi,\mathrm{bl}}})_{\overline{\mathbb F_q}} \to (\overline{{M^{\psi}}})_{\overline{\mathbb F_q}}$ is the induced map on the geometric special fibers. Since $b_{\overline s}$ is an isomorphism away from its fiber $Z_1$ over $y^{\mathrm{sg}}$, we deduce that
\begin{equation*}
(\mathrm R^i \Psi_{\eta}\Lambda)_{y} \simeq 
\begin{cases}
\Lambda & \text{if } i = 0 \text{ and } y \not = y^{\mathrm{sg}},\\
0 & \text{if } i>0 \text{ and } y \not = y^{\mathrm{sg}},\\
\mathrm H^i(Z_1,i_1^*\mathrm R \Psi_{\eta}^{\mathrm{bl}}\Lambda) & \text{if } y = y^{\mathrm{sg}}.
\end{cases}
\end{equation*}
Thus it remains to compute the cohomology of $Z_1$ with coefficients in $i_1^*\mathrm R \Psi_{\eta}^{\mathrm{bl}}\Lambda$.\\
Since $\mathrm R^i \Psi_{\eta}\Lambda$ has only two non-zero cohomology sheaves, the truncation functors yield a distinguished triangle
\begin{equation}\label{DistinguishedTriangle}
\Lambda[0] \to \mathrm R \Psi_{\eta}^{\mathrm{bl}}\Lambda \to i_{Q*}\Lambda(-1)[-1] \rightsquigarrow
\end{equation}
in $D_c^b((\overline{M^{\psi,\mathrm{bl}}})_{\overline{\mathbb F_q}}, \Lambda)$. Now, for every $K \in D_c^b((\overline{M^{\psi,\mathrm{bl}}})_{\overline{\mathbb F_q}}, \Lambda)$, there is a distinguished triangle 
\begin{equation}\label{GysinDistinguishedTriangle}
i_{1*}i_1^!K \to K \to \mathrm Rj_{1*}j_1^*K \rightsquigarrow,
\end{equation}
which is functorial in $K$. Applying this with $K = \Lambda[0], \mathrm R \Psi_{\eta}^{\mathrm{bl}}\Lambda$ and $i_{Q*}\Lambda(-1)[-1]$ successively and using functoriality with respect to the distinguished triangle \eqref{DistinguishedTriangle}, we obtain a large commutative diagram all of whose rows and columns are exact, see Figure \ref{figure1}. The columns are the Gysin sequences, ie. the long exact sequences induced by \eqref{GysinDistinguishedTriangle}. The rows are the long exact sequences induced by \eqref{DistinguishedTriangle}, and by the triangles obtained after applying $i_{1*}i_1^!$ and $\mathrm Rj_{1*}j_1^*$. In this Figure, we took into account the results of Section \ref{Section2} to replace vanishing cohomology groups with $0$. We also used the fact that $j_1^* \mathrm R \Psi_{\eta}^{\mathrm{bl}}\Lambda \simeq \Lambda[0]$, which holds since $U_1$ is in the smooth locus of $\overline{M^{\psi,\mathrm{bl}}}$. We single out the morphisms 
\begin{align*}
\alpha_i& :H^{i}(\overline{M^{\psi,\mathrm{bl}}},\Lambda) \to \mathrm H^{i}(\overline{M^{\psi,\mathrm{bl}}},\mathrm R \Psi_{\eta}^{\mathrm{bl}}\Lambda), & \beta_i: H^{i}(Z_1,i_1^!\Lambda) \to \mathrm H^{i}(Z_1,i_1^!\mathrm R \Psi_{\eta}^{\mathrm{bl}}\Lambda), \\
\varphi_{2i} & : \mathrm H^{2i}(\overline{M^{\psi,\mathrm{bl}}},\mathrm R \Psi_{\eta}^{\mathrm{bl}}\Lambda) \to \mathrm H^{2i}(U_1,\Lambda), & 
\end{align*}
as in the diagram.

\begin{figure}[h]
\begin{center}
\begin{tikzcd}
& & 0 \arrow[d] & 0 \arrow[d] & \\
0 \arrow[r] & \mathrm H^{2i-2}(Q)(-1) \arrow[r] \arrow[d,equal] & \mathrm H^{2i}(Z_1,i_1^!\Lambda) \arrow[r,"\beta_{2i}"] \arrow[d] & \mathrm H^{2i}(Z_1,i_1^!\mathrm R \Psi_{\eta}^{\mathrm{bl}}\Lambda) \arrow[r] \arrow[d] & 0 \\
0 \arrow[r] & \mathrm H^{2i-2}(Q)(-1) \arrow[r] & \mathrm H^{2i}(\overline{M^{\psi,\mathrm{bl}}},\Lambda) \arrow[r,"\alpha_{2i}"] \arrow[d,"j_1^*"] & \mathrm H^{2i}(\overline{M^{\psi,\mathrm{bl}}},\mathrm R \Psi_{\eta}^{\mathrm{bl}}\Lambda) \arrow[r] \arrow[d,"\varphi_{2i}"] & 0 \\
& 0 \arrow[r] & \mathrm H^{2i}(U_1,\Lambda) \arrow[r,equal] \arrow[d] & \mathrm H^{2i}(U_1,\Lambda) \arrow[r] \arrow[d] & 0 \\
& 0 \arrow[r] & \mathrm H^{2i+1}(Z_1,i_1^!\Lambda) \arrow[r,"\beta_{2i+1}"] \arrow[d] & \mathrm H^{2i+1}(Z_1,i_1^!\mathrm R \Psi_{\eta}^{\mathrm{bl}}\Lambda) \arrow[r,"\gamma"] \arrow[d] & \mathrm H^{2i}(Q,\Lambda)(-1) \arrow[d,equal] \\
& 0 \arrow[r] & \mathrm H^{2i+1}(\overline{M^{\psi,\mathrm{bl}}},\Lambda) \arrow[r,"\alpha_{2i+1}"] \arrow[d] & \mathrm H^{2i+1}(\overline{M^{\psi,\mathrm{bl}}},\mathrm R \Psi_{\eta}^{\mathrm{bl}}\Lambda) \arrow[r] \arrow[d] & \mathrm H^{2i}(Q,\Lambda)(-1) \\
& & 0 & 0 &
\end{tikzcd}
\end{center}
\caption{The commutative diagram with exact rows and columns.}\label{figure1}
\end{figure}

\begin{lemma}\label{CohomologyNearbyCycles}
We have 
\begin{equation*}
\mathrm H^{i}(\overline{M^{\psi,\mathrm{bl}}},\mathrm R \Psi_{\eta}^{\mathrm{bl}}\Lambda) = \begin{cases}
\Lambda(-\frac{i}{2}) & \text{if } 0 \leq i \leq 2(n-1) \text{ is even},\\
0 & \text{else}.
\end{cases}
\end{equation*}
\end{lemma}

\begin{proof}
Since $M^{\psi,\mathrm{bl}}$ is proper, we have
\begin{equation*}
\mathrm H^{i}(\overline{M^{\psi,\mathrm{bl}}},\mathrm R \Psi_{\eta}^{\mathrm{bl}}\Lambda) \simeq \mathrm H^i(M^{\psi,\mathrm{bl}}_{\overline{\eta}},\Lambda),
\end{equation*}
where $M^{\psi,\mathrm{bl}}_{\overline{\eta}}$ is the geometric generic fiber of $M^{\psi,\mathrm{bl}}$. It is isomorphic to the generic fiber of $M^{\psi}$, which itself is isomorphic to $\mathbb P^{n-1}$ according to \cite{pappas} Proof of Proposition 3.8. The result follows.
\end{proof}

In particular, the terms $\mathrm H^{2i+1}(\overline{M^{\psi,\mathrm{bl}}},\mathrm R \Psi_{\eta}^{\mathrm{bl}}\Lambda)$ and $\mathrm H^{2i+1}(\overline{M^{\psi,\mathrm{bl}}},\Lambda)$ in the bottom row of Figure \ref{figure1} vanish. It follows also that the map labelled $\gamma$ in the diagram is actually the zero map, and that $\beta_{2i+1}$ is an isomorphism. 

\begin{proposition}\label{Imagej_1}
The image of the composition 
\begin{equation*}
(\mathrm{pr}_{|U_1}^*)^{-1} \circ j_1^*: \mathrm H^{2i}(\overline{M^{\psi,\mathrm{bl}}},\Lambda) \to \mathrm H^{2i}(U_1,\Lambda) \xrightarrow{\sim} \mathrm H^{2i}(Q,\Lambda)
\end{equation*}
is included in $\iota_1^*(\mathrm H^{2i}(Z_1,\Lambda)) \subset \mathrm H^{2i}(Q,\Lambda)$. 
\end{proposition}

Unless $n$ is even and $i = \frac{n-2}{2}$, we actually have $\mathrm H^{2i}(Q,\Lambda) = \iota_1^*(\mathrm H^{2i}(Z_1,\Lambda))$ so that the statement is trivially true. Nonetheless, we give the proof for every $i$ as it is completely formal.

\begin{proof}
The open immersion $j_1$ is the composition of the embeddings 
\begin{equation*}
j_1 = i_2 \circ j_1': U_1 \hookrightarrow Z_1 \hookrightarrow \overline{M^{\psi,\mathrm{bl}}}.
\end{equation*}
Thus, by functoriality the restriction $j_1^*$ is given by the composition 
\begin{equation}\label{j_1}
j_1^* = (j_1')^* \circ i_2^*: \mathrm H^{2i}(\overline{M^{\psi,\mathrm{bl}}},\Lambda) \to \mathrm H^{2i}(Z_2,\Lambda) \to \mathrm H^{2i}(U_1,\Lambda).
\end{equation}
By definition, for $x \in \mathrm H^{2i}(\overline{M^{\psi,\mathrm{bl}}},\Lambda)$ the element $z := (\mathrm{pr}_{|U_1}^*)^{-1} \circ j_1^*(x)$ is characterized as the unique element of $\mathrm H^{2i}(Q,\Lambda)$ such that $(j_1')^*\circ \mathrm{pr}^*(z) = j_1^*(x)$. It follows that $i_2^*(x) - \mathrm{pr}^*(z)$ lies in $\mathrm{Ker}((j_1')^*)$. Now, consider the Gysin sequence associated to $\iota_2: Q \hookrightarrow Z_2 \hookleftarrow U_1:j_1'$. 
\begin{equation}\label{j_1'}
0 \to \mathrm H^{2i-2}(Q,\Lambda)(-1) \xrightarrow{\iota_{2*}} \mathrm H^{2i}(Z_2,\Lambda) \xrightarrow{(j_1')^*} \mathrm H^{2i}(U_1,\Lambda) \to 0.
\end{equation}
It follows that $\mathrm{Ker}((j_1')^*) = \mathrm{Im}(\iota_{2*})$. Thus, there exists some $y \in \mathrm H^{2i-2}(Q,\Lambda)(-1)$ such that $i_2^*(x) = \mathrm{pr}^*(z) + \iota_{2*}(y)$. Applying $\iota_2^*$, we obtain 
\begin{equation*}
i_Q^*(x) = z + \iota_2^*\circ \iota_{2*}(y) \in \mathrm H^{2i}(Q,\Lambda).
\end{equation*}
Indeed, we have $\mathrm{pr}\circ\iota_2 = \mathrm{id}_Q$, so that by functoriality we have $\iota_2^*\circ \mathrm{pr}^* = \mathrm{id}$ on $\mathrm H^{2i}(Q,\Lambda)$. Since $i_Q^{*}(x) = \iota_1^*\circ i_{1}^*(x)$ and $\iota_2^*\circ \iota_{2*}(y) = - \iota_1^*\circ \iota_{1*}(y)$, it follows that $z \in \iota_1^*(\mathrm H^{2i}(Z_1,\Lambda))$.
\end{proof}

\begin{remark}
The fact that $\iota_2^*\circ \iota_{2*} + \iota_1^*\circ \iota_{1*} = 0$ can be seen by noticing that $\iota_i^* \circ \iota_{i*}$ is equal to the cup product with $c_1(\mathcal N_{Q/Z_i})$, the first Chern class of the normal bundle of $Q$ in $Z_i$, cf. Theorem 4.1 of \cite{SGA5} (exposé VII of SGA5). Since $\overline{M^{\psi,\mathrm{bl}}}$ is a principal divisor in $M^{\psi,\mathrm{bl}}$, one may check that $c_1(\mathcal N_{Z_1}(Q)) = - c_1(\mathcal N_{Z_2}(Q))$.
\end{remark}

\begin{proposition}\label{SurjectivePhi}
The morphism $\varphi_{2i}$ is surjective for all $0 \leq i \leq n-1$, except when $n$ is even and $i = \frac{n-2}{2}$.
\end{proposition}

\begin{proof}
By the diagram, we have $\varphi_{2i} \circ \alpha_{2i} = j_1^*:\mathrm H^{2i}(\overline{M^{\psi,\mathrm{bl}}},\Lambda) \to \mathrm H^{2i}(U_1,\Lambda)$, thus $\mathrm{Im}(j_1^*) \subset \mathrm{Im}(\varphi_{2i})$. Therefore it is enough to prove that $j_1^*$ is surjective. Recall from \eqref{j_1} that $j_1^* = (j_1')^* \circ i_2^*$, and from \eqref{j_1'} that $(j_1')^*$ is always surjective. Let us now consider the dual Gysin sequence associated to $i_2:Z_2 \hookrightarrow \overline{M^{\psi,\mathrm{bl}}} \hookleftarrow U_2:j_2$, ie. the following long exact sequence.
\begin{equation}\label{i_2}
\ldots \to \mathrm H_c^{2i}(U_2,\Lambda) \to \mathrm H^{2i}(\overline{M^{\psi,\mathrm{bl}}},\Lambda) \xrightarrow{i_2^*} \mathrm H^{2i}(Z_2,\Lambda) \to \mathrm H_c^{2i+1}(U_2,\Lambda) \to \ldots
\end{equation}
By Poincaré duality, it follows from Lemma \ref{CohomologyU2} that $\mathrm H_c^i(U_2,\Lambda) = 0$ unless when $i = 2(n-1)$, and when $n$ is even and $i = n-1$. Thus, the restriction $i_2^*:\mathrm H^{2i}(\overline{M^{\psi,\mathrm{bl}}},\Lambda) \to \mathrm H^{2i}(Z_2,\Lambda)$ is an isomorphism for all $0 \leq i < n-1$ with $i \not = \frac{n-2}{2}$ when $n$ is even, and it is surjective when $i = n-1$. By composition, it follows that $j_1^*:\mathrm H^{2i}(\overline{M^{\psi,\mathrm{bl}}},\Lambda) \to \mathrm H^{2i}(U_1,\Lambda)$ is surjective for all $0 \leq i \leq n-1$ with $i \not = \frac{n-2}{2}$ when $n$ is even.\end{proof} 

%Recall that we have $\mathrm{pr}\circ \iota_2 = \mathrm{id}_Q$. By functoriality, it follows that the compositions
%\begin{align*}
%\iota_2^* \circ \mathrm{pr}^*: \mathrm H^{2i}(Q,\Lambda) \to \mathrm H^{2i}(Z_2,\Lambda) \to \mathrm H^{2i}(Q,\Lambda), & & \mathrm{pr}_* \circ \iota_{2*} : \mathrm H^{2i-2}(Q,\Lambda)(-1) \to \mathrm H^{2i}(Z_2,\Lambda) \to \mathrm H^{2i-2}(Q,\Lambda)(-1),
%\end{align*}
%are the identity respectively on $\mathrm H^{2i}(Q,\Lambda)$ and on $\mathrm H^{2i-2}(Q,\Lambda)(-1)$. For $x$ in the former group and $y$ in the latter, let us write $z := (\mathrm{pr}^*(x),\zeta \smile \mathrm{pr}^*(y)) \in \mathrm H^{2i}(Z_2,\Lambda)$. Then we have $\mathrm{pr}_*(z) = y$ and $\iota_2^*(z) = x$. Since the composition $\iota_2^*\circ \iota_{2*}$ corresponds to the cup product with $c_1(\mathcal N) \in \mathrm H^{2}(Q,\Lambda)(-1)$. It follows that $\iota_{2*}(y) = (\mathrm{cl}_{Z_2}(Q)\smile y, \zeta \smile y)$ for 

% In this case, we have 
%\begin{equation*}
%\mathrm H^{n-2}(U_1,\Lambda) \simeq \mathrm H^{n-2}(Q,\Lambda) \simeq \iota_1^*(\mathrm H^{n-2}(Z_1,\Lambda)) \oplus \mathrm H^{n-2}_{\mathrm{prim}}(Q,\Lambda),
%\end{equation*} 
%and the image of $\varphi_{n-2}$ is identified with $\iota_1^*(\mathrm H^{n-2}(Z_1,\Lambda))$.

\begin{corollary}\label{VanishingRPsi}
We have $\mathrm H^0(Z_1,i_1^*\mathrm R \Psi_{\eta}^{\mathrm{bl}}\Lambda) = \Lambda$ and $\mathrm H^i(Z_1,i_1^*\mathrm R \Psi_{\eta}^{\mathrm{bl}}\Lambda) = 0$ for all $i \geq 1$ such that $i \not \in \{n-1,n\}$ when $n$ is even.
\end{corollary}

\begin{proof}
Recall the duality functor $D$ introduced in the Notations \ref{Notations}. By Poincaré duality, we have 
\begin{equation*}
\mathrm H^i(Z_1,i_1^*\mathrm R \Psi_{\eta}^{\mathrm{bl}}\Lambda) = \mathrm H^{-i}(Z_1,Di_1^*\mathrm R \Psi_{\eta}^{\mathrm{bl}}\Lambda)^{\vee},
\end{equation*}
where $(\cdot)^{\vee}$ denotes the $\Lambda$-linear dual. Moreover, we have 
\begin{equation*}
Di_1^*\mathrm R \Psi_{\eta}^{\mathrm{bl}} \Lambda = i_1^!D\mathrm R \Psi_{\eta}^{\mathrm{bl}} \Lambda = i_1^!\mathrm R \Psi_{\eta}^{\mathrm{bl}} D \Lambda = i_1^!\mathrm R \Psi_{\eta}^{\mathrm{bl}} \Lambda (n-1)[2(n-1)].
\end{equation*}
Indeed, the nearby cycles $\mathrm R\Psi$ commute with duality by \cite{HansenScholze} Corollary 4.2. The last equality follows from the smoothness of the generic fiber of $M^{\psi,\mathrm{bl}}$. It follows that 
\begin{equation}\label{PoincaréDuality}
\mathrm H^i(Z_1,i_1^*\mathrm R \Psi_{\eta}^{\mathrm{bl}}\Lambda) = \mathrm H^{-i+2(n-1)}(Z_1,i_1^!\mathrm R \Psi_{\eta}^{\mathrm{bl}}\Lambda)^{\vee}(1-n).
\end{equation}
Let $0 \leq i \leq n-1$, and assume that $i \not = \frac{n-2}{2}$ when $n$ is even. By Proposition \ref{SurjectivePhi}, the morphism $\varphi_{2i}$ is surjective. From the diagram in Figure \ref{figure1}, it follows that $\mathrm H^{2i+1}(Z_1,i_1^!\mathrm R \Psi_{\eta}^{\mathrm{bl}}\Lambda) = 0$.\\
Assume $i \not = n-1$. Ignoring the Tate twists, we have $\mathrm H^{2i}(U_1,\Lambda) \simeq \Lambda$ by Lemmas \ref{CohomologyOfQ} and \ref{CohomologyU1}, and $\mathrm H^{2i}(\overline{M^{\psi,\mathrm{bl}}},\mathrm R \Psi_{\eta}^{\mathrm{bl}}\Lambda) \simeq \Lambda$ by Lemma \ref{CohomologyNearbyCycles}. Since any surjective map $\Lambda \to \Lambda$ of $\Lambda$-modules is an automorphism, $\varphi_{2i}$ is also injective. Thus $\mathrm H^{2i}(Z_1,i_1^!\mathrm R \Psi_{\eta}^{\mathrm{bl}}\Lambda) = \mathrm{Ker}(\varphi_{2i}) = 0$. \\
\sloppy Eventually, if $i = n-1$ we have $\mathrm H^{2(n-1)}(U_1,\Lambda) = 0$ so that $\mathrm H^{2(n-1)}(Z_1,i_1^!\mathrm R \Psi_{\eta}^{\mathrm{bl}}\Lambda) \simeq H^{2(n-1)}(\overline{M^{\psi,\mathrm{bl}}},\mathrm R \Psi_{\eta}^{\mathrm{bl}}\Lambda) \simeq \Lambda(1-n)$.
\end{proof}

The proof of Theorem \ref{ComputationOfNearbyCycles} is now over when $n$ is odd, and all that remains is to determine $\mathrm H^{n-1}(Z_1,i_1^*\mathrm R \Psi_{\eta}^{\mathrm{bl}}\Lambda)$ and $\mathrm H^{n}(Z_1,i_1^*\mathrm R \Psi_{\eta}^{\mathrm{bl}}\Lambda)$ when $n$ is even. 

\begin{proposition}
Assume that $n$ is even. Then $\mathrm H^{n}(Z_1,i_1^*\mathrm R \Psi_{\eta}^{\mathrm{bl}}\Lambda) = 0$ and $\mathrm H^{n-1}(Z_1,i_1^*\mathrm R \Psi_{\eta}^{\mathrm{bl}}\Lambda) \simeq \Lambda[\epsilon q^{\frac{n}{2}}]$.
\end{proposition}

\begin{proof}
The sheaf $\mathrm R \Psi_{\eta}\Lambda[n-1]$ is an object of $\mathrm{Perv}((\overline{M^{\psi}})_{\overline{\mathbb F_q}},\Lambda)$, since $\Lambda[n-1]$ is perverse on the generic fiber and the nearby cycles preserve perversity by \cite{HansenScholze} Lemma 6.3. By definition of perversity and since $y^{\mathrm{sg}}$ is closed in $\overline{M^{\psi}}$, we know that $\mathrm H^{k+n-1}i_{\overline{y^{\mathrm{sg}}}}^*\mathrm R \Psi_{\eta}\Lambda = 0$ for all $k \geq 1$. In particular, we have
\begin{equation*}
\mathrm H^{n}(Z_1,i_1^*\mathrm R \Psi_{\eta}^{\mathrm{bl}}\Lambda) \simeq (\mathrm R^{n} \Psi_{\eta}\Lambda)_{\overline{y^{\mathrm{sg}}}} = 0.
\end{equation*}

By Poincaré duality \eqref{PoincaréDuality}, it follows that $\mathrm H^{n-2}(Z_1,i_1^!\mathrm R \Psi_{\eta}^{\mathrm{bl}}\Lambda) = 0$ when $n$ is even, so that $\varphi_{n-2}$ is injective. The codomain of $\varphi_{n-2}$ is 
\begin{equation*}
\mathrm H^{n-2}(U_1,\Lambda) \simeq \mathrm H^{n-2}(Q,\Lambda) \simeq \iota_1^*(\mathrm H^{n-2}(Z_1,\Lambda))\oplus \mathrm H^{n-2}_{\mathrm{prim}}(Q,\Lambda),
\end{equation*}
and in Proposition \ref{Imagej_1} we proved that $\mathrm{Im}(\varphi_{n-2}) = \mathrm{Im}(j_1^*)$ is included in $\iota_1^*(\mathrm H^{n-2}(Z_1,\Lambda))$. It follows that 
\begin{equation}\label{TheTorsionPart}
\mathrm H^{n-1}(Z_1,i_1^!\mathrm R \Psi_{\eta}^{\mathrm{bl}}\Lambda) \simeq \mathrm H^{n-2}_{\mathrm{prim}}(Q,\Lambda) \oplus \frac{\iota_1^*(\mathrm H^{n-2}(Z_1,\Lambda))}{\mathrm{Im}(\varphi_{n-2})}.
\end{equation}
Forgetting about the Tate twists, we have $\iota_1^*(\mathrm H^{n-2}(Z_1,\Lambda)) \simeq \Lambda$ and $\mathrm{Im}(\varphi_{n-2}) \simeq \mathrm H^{n-2}(\overline{M^{\psi,\mathrm{bl}}},\mathrm R \Psi_{\eta}^{\mathrm{bl}}\Lambda) \simeq \Lambda$. Thus, the quotient module in the right-hand side of \eqref{TheTorsionPart} is isomorphic to $\Lambda/a\Lambda$ for some non-zero divisor $a \in \Lambda$. Therefore, it is killed when taking the $\Lambda$-linear dual. By Poincaré duality \eqref{PoincaréDuality} we have
\begin{equation*}
\mathrm H^{n-1}(Z_1,i_1^*\mathrm R \Psi_{\eta}^{\mathrm{bl}}\Lambda) \simeq \mathrm H^{n-2}_{\mathrm{prim}}(Q,\Lambda)^{\vee}(1-n) \simeq \Lambda[\epsilon q^{\frac{n}{2}}].
\end{equation*}
\end{proof}

This concludes the proof of Theorem \ref{ComputationOfNearbyCycles}.

\begin{remark}\label{LastRemark}
As mentioned in the Introduction, Theorem \ref{ComputationOfNearbyCycles} allows us to recover Krämer's computation of the semisimple trace of the Frobenius. Indeed, let $x \in M^{\psi}(\mathbb F_q)$ and let $\overline x$ be the geometric point above it. We define 
\begin{equation*}
\mathrm{Tr}(\mathrm{Frob},(\mathrm R\Psi_{\eta}\Lambda)_{\overline x}) := \sum_{k \geq 0} (-1)^k \mathrm{Tr}(\mathrm{Frob},(\mathrm R^k\Psi_{\eta}\Lambda)_{\overline x}).
\end{equation*}
By our computations, we have 
\begin{equation*}
\mathrm{Tr}(\mathrm{Frob},(\mathrm R\Psi_{\eta}\Lambda)_{\overline x}) = \begin{cases}
1 & \text{if } x \not = y^{\mathrm{sg}} \text{ or if } n \text{ is odd},\\
1 - \epsilon q^{\frac{n}{2}} & \text{if } x = y^{\mathrm{sg}} \text{ and } n \text{ is even}.
\end{cases}
\end{equation*}
In \cite{kramer}, the author only consider the alternating form $\psi_0$ given by the matrix $\mathrm{Diag}(1,\ldots,1)$ in a certain basis. The resulting hermitian space $(V,(\cdot,\cdot))$ is split if and only if $(-1)^{\frac{n(n-1)}{2}} \in \mathrm{Norm}_{E/E_0}(E^{\times})$, which is equivalent to $\left(\frac{-1}{q}\right)^{\frac{n(n-1)}{2}} = 1$. Now consider the space $\mathbb F_q^n$ with the symmetric bilinear pairing $x \mapsto \sum x_i^2$. Assume that $n$ is even. Then $\mathbb F_q^{n}$ is hyperbolic (ie. isomorphic to a direct sum of $\frac{n}{2}$-copies of the hyperbolic plane) if and only if $\left(\frac{-1}{q}\right)^{\frac{n(n-1)}{2}} = 1$. Putting things together, we recover the formula
\begin{equation*}
\mathrm{Tr}(\mathrm{Frob},(\mathrm R\Psi_{\eta}\Lambda)_{\overline y^{\mathrm{sg}}}) = \begin{cases}
1 & \text{if } n \text{ is odd},\\
1 - q^{\frac{n}{2}} & \text{if } n \text{ is even and } \mathbb F_q^n \text{ is hyperbolic},\\
1 + q^{\frac{n}{2}} & \text{if } n \text{ is even and } \mathbb F_q^n \text{ is not hyperbolic},
\end{cases}
\end{equation*}
which is valid for $\psi = \psi_0$.
\end{remark}

\bibliographystyle{amsplain}
\bibliography{biblio.bib}

\end{document}